\documentclass[a4paper,oneside,12pt]{article}
\usepackage[latin1]{inputenc}
\usepackage{lmodern}
\usepackage[T1]{fontenc}
\usepackage{textcomp}
\usepackage[french,english]{babel}
\usepackage[a4paper,vmargin={2.5cm,3.5cm},hmargin=2cm]{geometry}
\usepackage[pdftex]{hyperref}
\usepackage[expansion=false]{microtype}
\usepackage[pdftex]{color,graphicx}

\usepackage{amsmath,amsfonts,amssymb,amsthm,mathrsfs}

\theoremstyle{plain}
\newtheorem{theorem}{Theorem}
\newtheorem{corollary}{Corollary}
\newtheorem{proposition}{Proposition}
\newtheorem{lemma}{Lemma}
\theoremstyle{definition}
\newtheorem{definition}{Definition}
\newtheorem{step}{Step}
\theoremstyle{remark}
\newtheorem*{remark}{Remark}

\hypersetup{
    pdfauthor   = {Laurent Menard},%
    pdftitle    = {The two uniform infinite quadrangulations of the plane have the same law}}

\author{Laurent {\sc M\'enard}\footnote{\texttt{laurent.menard@normalesup.org}}\\
        D\'epartement de Math\'ematiques\\
        Universit\'e Paris-Sud\\
        91405 Orsay Cedex, France}
\title{\bf THE TWO UNIFORM INFINITE QUADRANGULATIONS OF THE PLANE HAVE THE SAME LAW}
\date{{\small December 2008}}

\begin{document}

\maketitle

\selectlanguage{french}

\begin{abstract}
On d\'emontre que les quadrangulations al\'eatoires infinies
uniformes d\'efinies respectivement par Chassaing-Durhuus et par
Krikun ont la m\^eme loi.
\end{abstract}

\selectlanguage{english}

\begin{abstract}
We prove that the uniform infinite random quadrangulations defined
respectively by Chassaing-Durhuus and Krikun have the same distribution.
\end{abstract}

\noindent {\bf AMS classification:} 60C05, 60J80, 05C30.

\section{Introduction}

Planar maps are proper embeddings of connected graphs in the two-dimensional sphere $\mathbb{S}^2$. Their combinatorial properties
have been studied by Tutte \cite{Tutte}
and many others. Planar maps have recently drawn much attention in the theoretical physics literature as models of random surfaces, especially in the setting of the theory of two-dimensional quantum gravity (see in particular the book \cite{2DQG}). A powerful tool to study these objects is the encoding of planar maps in terms of labelled trees, which was first introduced by Cori and Vauquelin in \cite{CV} and was much developed in Schaeffer's thesis \cite{Schaeffer} (see also Bouttier, Di Francesco and Guitter \cite{BDG} for a generalized version of this encoding). This correspondence between planar maps and trees makes it possible to derive certain asymptotics of large random planar maps in terms of continuous random trees (see the work of Chassaing and Schaeffer \cite{CS}) and to define a Brownian map (see Marckert and Mokkadem \cite{MaMo}) which is a continuous random metric space conjectured to be the scaling limit of various classes of planar maps (see the papers by Marckert and Miermont \cite{MaMi}, Le Gall \cite{LG}, Le Gall and Paulin \cite{LGP}). This approach has led to new asymptotic properties of large planar maps.

\bigskip

Another point of view is to study properties of random infinite planar maps, more precisely to study probability measures on certain classes of infinite planar maps, which are uniform in some sense. This has been done by Angel and Schramm \cite{AS} who introduced a uniform infinite triangulation of the plane, later studied by Angel \cite{A1,A2} and Krikun \cite{Krikun}.

In the present paper, we are interested in infinite random planar quadrangulations. Recall that Schaeffer's bijection (see e.g. \cite{CS}) yields a one-to-one correspondence between rooted planar quadrangulations with $n$ faces and well-labelled trees with $n$ edges. Then, there are two natural ways to define a uniform infinite quadrangulation
of the sphere: one as the local limit of uniform finite quadrangulations as their size goes to infinity, and one going through Schaeffer's bijection and using local limits of uniform well-labelled trees. The first approach is developed in Krikun \cite{Kr} while the second one is developed in Chassaing and Durhuus \cite{CD}. The topologies in which the uniform finite quadrangulations converge to the infinite object differ in the two cases: in the Chassaing-Durhuus paper, the topology on quadrangulations is induced by Schaeffer's bijection and the natural topology of local convergence of rooted trees, while the topology used in Krikun's paper is the natural topology of local convergence of rooted planar maps. Therefore, the two uniform infinite random quadrangulations defined in these papers are a priori two different objects. The goal of the paper is to show that these two definitions coincide. This result is stated in Theorem \ref{equality} below. Note that our work also gives an alternative approach to
Theorem 1 of Krikun \cite{Kr}: independently of the results of \cite{Kr}, Theorem \ref{equality} shows that
the uniform probability measure on the space of all rooted planar quadrangulations with $n$ faces
converges as $n\to\infty$ to a probability measure on the space of infinite quadrangulations,
in the sense of the metric used in \cite{Kr}. 
\bigskip

Let us briefly explain the main point of our argument. Consider a sequence of (deterministic or random) finite well-labelled trees $\theta_n$, that converges as $n \to \infty$ towards an infinite well-labelled tree $\theta_{\infty}$, in the sense that, for every $k \geqslant 1$, the restriction of $\theta_n$ to the 
first $k$ generations is equal to the same restriction of $\theta_{\infty}$,
when $n$ is sufficiently large. Let $Q_n$ be the quadrangulation associated with $\theta_n$ via Schaeffer's bijection and let $Q_{\infty}$ be the infinite quadrangulation associated with $\theta_{\infty}$ via the extension of Schaeffer's bijection that is presented in Subsection \ref{subsec:schaeffer} below ($\theta_{\infty}$ needs to satisfy certain properties so that this makes sense). Then it is not always true that $Q_{\infty}$ is the local limit of $Q_n$ as $n \to \infty$. The problem comes from the fact that $\theta_n$ may have small labels at generations larger than $k(n)$ with $k(n) \to \infty$. Note that this problem may occur even if one knows that $\theta_{\infty}$ has finitely many labels smaller than $K$, for every integer $K$ (the latter property holds for the uniform infinite well-labelled tree thanks to the estimates of \cite{CD}, see Proposition \ref{etiquettes} below). Nonetheless, in the case when $\theta_n$ is uniformly distributed over all well-labelled trees with $n$ edges, the preceding phenomenon does not occur: for every fixed $R>0$, the probability that $\theta_n$ has a label less than $R$ above generation $S$ tends to $0$ as $S \to \infty$, uniformly in $n$. This uniform estimate is stated in Proposition \ref{nseuil} below.

We can combine this estimate with the following combinatorial argument. If two well-labelled trees coincide up to generation $S$, then the associated quadrangulations are also the same within distance $R$ from the root, where $R$ is essentially the minimum label above generation $S$ in either tree. See Proposition \ref{facegen} below for a precise statement.

\bigskip

The paper is organized as follows: Section 2 gives some notation and an
extension of Schaeffer's bijection to the infinite case; Section 3
presents the two different definitions of the uniform infinite quadrangulation; and
Section 4 contains the key estimates that allow us to prove that these definitions actually lead to the same object.

\section{Preliminaries}
\label{sec:prel}

\subsection{Spatial trees}
\label{subsec:spatial trees}

Throughout this work we will use the standard formalism on planar trees as found in \cite{Neveu}. Let \[ \mathcal{U} = \bigcup_{n=0}^{\infty} \mathbb{N}^n\] where by convention $\mathbb{N} = \{ 1,2, \ldots \}$ and $\mathbb{N}^0 = \left\{ \emptyset \right\}$. An element $u$ of $\mathcal{U}$ is thus a finite sequence of positive integers. If $u, v \in \mathcal{U}$, $uv$ denotes the concatenation of $u$ and $v$. If $v$ is of the form $uj$ with $j \in \mathbb{N}$, we say that $u$ is the \emph{parent} of $v$ or that $v$ is a \emph{child} of $u$. More generally, if $v$ is of the form $uw$ for $u,w \in \mathcal{U}$, we say that $u$ is an \emph{ancestor} of $v$ or that $v$ is a \emph{descendant} of $u$. A \emph{rooted planar tree} $\tau$ is a subset of $\mathcal{U}$ such that
\begin{enumerate}
\item $\emptyset \in \tau$  ( $\emptyset$ is called the \emph{root} of $\tau$),
\item if $v \in \tau$ and $v \neq \emptyset$, the parent of $v$ belongs to $\tau$
\item for every $u \in \mathcal{U}$ there exists $k_u(\tau) \geqslant 0$ such that $uj \in \tau$ if and only if $j \leqslant k_u(\tau)$.
\end{enumerate}
The edges of $\tau$ are the pairs $(u,v)$, where $u, v \in \tau$ and $u$ is the father of $v$. $|\tau|$ denotes the number of edges of $\tau$ and is called the size of $\tau$. $h(\tau)$ denotes the maximal generation of a vertex in $\tau$ and is called the height of $\tau$. We denote by $\mathbf{T}_n$ the set of all rooted planar trees of size $n$ and by $\mathbf{T}_{\infty}$ the set of all infinite rooted planar trees. Then $\mathbf{T} = \cup_{n=0}^{\infty} \mathbf{T}_n$ is the set of all  finite rooted planar trees and $\overline{\mathbf{T}} = \mathbf{T} \cup \mathbf{T}_{\infty}$ is the set of all rooted (finite or infinite) planar trees. A \emph{spine} of a tree $\tau$ is an infinite linear sub-tree of $\tau$ starting from its root.

\bigskip

A \emph{rooted labelled tree} (or spatial tree) is a pair $\theta = (\tau, (\ell(u))_{u \in \tau})$  that consists of a planar tree $\tau$ and a collection of integer labels assigned to the vertices of $\tau$ such that if $u,v \in \tau$ and $v$ is a child of $u$, then $|\ell(u) - \ell(v)| \leqslant 1$.
For every $l \in \mathbb{Z}$, we denote by $\overline{\mathbf{T}}^{(l)}$ the set of all spatial trees for which $\ell(\emptyset) = l$, by $\mathbf{T}_{\infty}^{(l)}$ the set of all such trees with an infinite number of edges, by $\mathbf{T}_n^{(l)}$ the set of all such trees with $n$ edges and by $\mathbf{T}^{(l)}$ the set of all such trees with finitely many vertices. Similarly as before, $\mathbf{T}^{(l)} = \bigcup_{n=0}^{\infty} \mathbf{T}_n^{(l)}$.

If $\ell(\emptyset) = l$ and in addition $\ell(u) \geqslant 1$ for every vertex $u$ of $\tau$, we say that $\theta$ is an $l$-well-labelled tree. The corresponding sets of spatial trees are denoted by $\overline{\mathbb{T}}^{(l)}$, $\mathbb{T}^{(l)}$, $\mathbb{T}_{\infty}^{(l)}$ and $\mathbb{T}_n^{(l)}$. For $l= 1$ we will simply say well-labelled tree and denote the corresponding sets by $\overline{\mathbb{T}}$, $\mathbb{T}$, $\mathbb{T}_{\infty}$ and $\mathbb{T}_n$.

\begin{figure}[!t]
\begin{center}
\includegraphics[width=\textwidth]{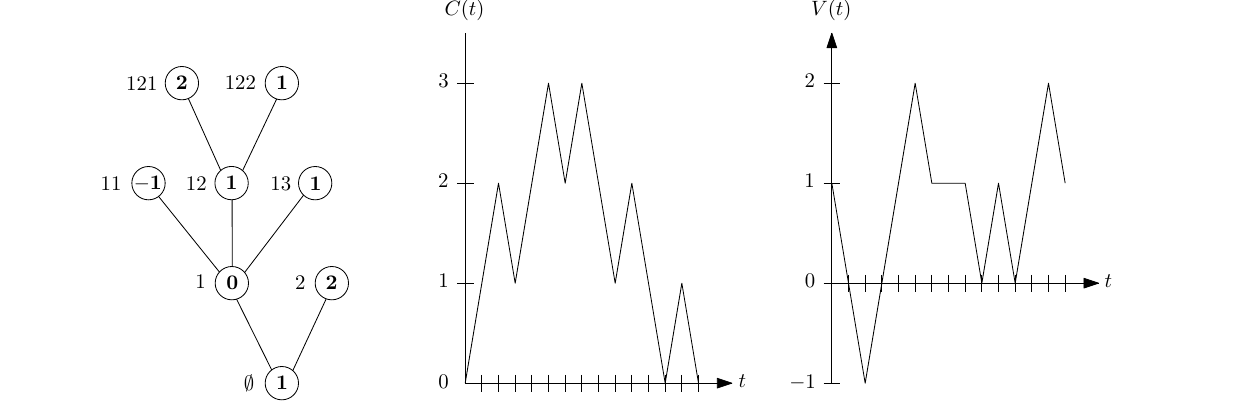}
\end{center}
\caption{A spatial tree and its pair of contour functions $(C,V)$.}
\label{fig:contour}
\end{figure}

A finite spatial tree $\omega = (\tau, \ell)$ can be coded by a pair $(C,V)$, where $C = (C(t))_{0 \leqslant t \leqslant 2|\tau|}$ is the contour function of $\tau$ and $V = (V(t))_{0 \leqslant t \leqslant 2|\tau|}$ is the spatial contour function of $\omega$ (see Figure \ref{fig:contour}). To define these contour functions, let us consider a particle which, starting from the root, traverses the tree along its edges at speed one. When leaving a vertex, the particle visits the first non visited child of this vertex if there is such a child, or returns to the parent of this vertex. Since all edges will be crossed twice, the total time needed to explore the tree is $2 |\tau|$. For every $t \in [0, 2|\tau|]$, $C(t)$ denotes the distance from the root of the position of the particle. In addition if $t \in [0, 2 |\tau|]$ is an integer, $V(t)$ denotes the label of the vertex that is visited at time $t$. We then complete the definition of $V$ by interpolating linearly between successive integers. See Figure \ref{fig:contour} for an example. A spatial tree is uniquely determined by its pair of contour functions.

\bigskip

To conclude this section, let us introduce some relevant notation. If $\omega = (\tau,\ell)$ is a labelled tree, $|\omega| = |\tau|$ is the size of $\omega$, $h(\omega) = h(\tau)$ is the height of $\omega$ and, for $S \geqslant 0$, $g_S(\omega)$ is the set of all vertices of $\omega$ at generation $S$. Finally, for every $l \in \mathbb{N}$, we let $N_l(\omega)$ denote the number of vertices of $\omega$ that have label $l$. We define $\mathscr{S}$ as the set of all trees of $\overline{\mathbb{T}}$ that have at most one spine, and for which labels takes each integer value a finite number of times:
\begin{equation}
\label{spine}
\mathscr{S} = \left\{ \omega \in \mathbb{T}_{\infty}: \, \forall l \geqslant 1, \, N_l (\omega) < \infty \text{ and $\omega$ has a unique spine}\right\} \cup \mathbb{T}.
\end{equation}

\subsection{Planar maps and quadrangulations}
\label{subsec:quadrangulations}

Consider a proper embedding of a finite connected graph in the sphere $\mathbb{S}^2$ (loops and multiple edges are allowed).  A (finite) \emph{planar map} is an equivalent class of such embedded graphs with respect to orientation preserving homeomorphisms of the sphere. A planar map is \emph{rooted} if it has a distinguished oriented edge, and the origin of the root is called the root vertex. In what follows, planar maps are always rooted even if this is not mentioned explicitly. The set of vertices will always be equipped with the graph distance. The faces of the map are the connected components of the complement of the union of its edges. A finite planar map is a \emph{quadrangulation} if all its faces have degree $4$.

For every integer $n \geqslant 1$ we let $\mathbf{Q}_n$ denote the set of all (rooted) quadrangulations with $n$ faces and $\mathbf{Q} = \bigcup_{n \geqslant 1} \mathbf{Q}_n$ denote the set of finite quadrangulations. Each set $\mathbf{Q}_n$ is in bijective correspondence with the set $\mathbb{T}_n$ by Schaeffer's bijection \cite{CV,Schaeffer}. There is no bijection between infinite well-labelled trees and infinite quadrangulations, but Schaeffer's correspondence has been extended to $\mathscr{S}$ in \cite{CD}. To discuss this extention, we first have to define precisely what we mean by an infinite quadrangulation. To this end we recall some definitions of \cite{AS,CD} in a slightly different form.

\bigskip
Throughout this work, we consider only infinite graphs such that the degree 
of every vertex is finite.
Consider a proper embedding of an infinite graph in the plane $\mathbb{R}^2$. We say that this embedding is locally finite if every compact subset of $\mathbb{R}^2$ intersects only finitely many edges.

\begin{definition}
\label{acceptable}
An infinite planar map $\mathcal{M}$ is an equivalent class of locally finite embeddings of an infinite graph in $\mathbb{R}^2$, with respect to orientation preserving homeomorphisms of the plane.
\end{definition}

The faces of an infinite planar map $\mathcal{M}$ are the bounded connected components of the complement of the union of its edges.
With this definition, every edge of $\mathcal{M}$ is not necessarily adjacent to a face; for example infinite trees have only one ''face'' of infinite degree, which is not a face in the sense of the previous definition. This motivates the next definition.

\begin{definition}
A \emph{regular} infinite planar map is  an infinite planar map such that every connected component of the complement of the union of its edges is bounded.
\end{definition}

In a \emph{regular} infinite planar map, every edge is either shared by two faces or appears twice in the border of a face.

\begin{remark}
With the previous definitions, an infinite tree can be embedded as an infinite planar map in $\mathbb{R}^2$, but not as a regular infinite planar map.
\end{remark}

\begin{definition}
\label{Qinf}
An \emph{infinite planar quadrangulation} is a regular infinite planar map having every face bordered by four-sided polygons. A rooted infinite quadrangulation is an infinite quadrangulation with a distinguished oriented edge $(v_0,v_1)$ called the root of the quadrangulation; $v_0$ is called the root vertex of the quadrangulation. We denote by $\overline{\mathbb{Q}}$ the set of all (finite or infinite) rooted planar quadrangulations and we have the self-evident decomposition $\overline{\mathbb{Q}} = \mathbf{Q}\cup \mathbf{Q}_{\infty}$.
\end{definition}

\subsection{Schaeffer's correspondence}
\label{subsec:schaeffer}

We are now going to describe the extension of Schaeffer's correspondence to the set $\mathscr{S}$. We refer to Section 6 of \cite{CD} for details and proofs.

With every infinite well-labelled tree $\omega \in \mathscr{S}$ we will associate an infinite planar 
quadrangulation $\Phi(\omega)$. We identify $\mathbb{S}^2$ with the set $\mathbb{R}^2 \cup \{\infty \}$, and we fix an infinite tree $\omega \in \mathscr{S}$.
We can also fix an embedding of $\omega$ into $\mathbb{R}^2$ as in Definition \ref{acceptable} above. 
We root $\omega$ at the edge between vertices $\emptyset$ and $1$.
Let $F_0$ denote the complement of 
the union of edges of $\omega$ in $\mathbb{S}^2$.

\begin{definition}
A \emph{corner} of $F_0$ is a sector between two consecutive edges around a vertex. The label of a corner is the label of the corresponding vertex.
\end{definition}

A vertex of degree $d$ defines $d$ corners and a tree $\omega \in \mathscr{S}$ has a finite number $C_k(\omega) \geqslant N_k(\omega)$ of corners with label $k$. The map $\Phi(\omega)$ is defined in three steps.

\begin{figure}[!t]
\begin{center}
\includegraphics[width=\textwidth]{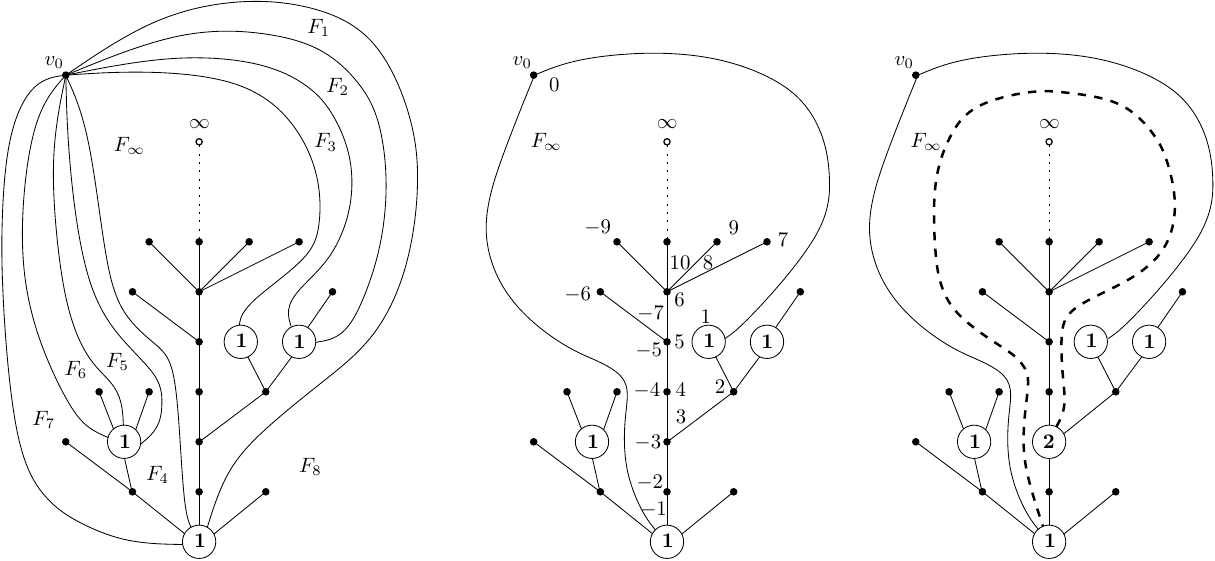}
\end{center}
\caption{Left: step 1, edges are added between $v_0$ and corners with label $1$. Middle: step 2, numbering of a few corners in $F_{\infty}$. Right: step 2, a chord between the two sides of the spine.}
\label{figschaeffer}
\end{figure}

\begin{step}[see Figure \ref{figschaeffer}]
A vertex $v_0$ with label $0$ is added in $F_0 \setminus \{\infty\}$ and one edge is added between this vertex and each of the $C_1(\omega)$ corners with label $1$. The new root is taken to be the edge that connects $v_0$ to the corner before the root edge of $\omega$.
\end{step}

\begin{remark}
Notice that the construction in step 1 is possible because $\omega$ has at most one spine.
\end{remark}

After step 1, a uniquely defined rooted infinite planar map $\mathcal{M}_0$ with $C_1(\omega) -1$ faces is obtained (in the sense of the definitions of Section \ref{subsec:quadrangulations}, in particular, the faces are bounded subsets of $\mathbb{R}^2$). Notice that each face of $\mathcal{M}_0$ has a unique corner with label $0$ and two corners with label $1$. Such a face is bordered by the two edges joining $v_0$ to the two corners with label $1$ and, in the case where the two corners with label $1$ correspond to two different vertices, the unique injective path in the tree between these two vertices with label $1$.

It is natural to consider the complement of $\mathcal{M}_0$ and its faces as an additional face of infinite degree. Let us denote this face by $F_{\infty}$. It possesses a unique corner with label $0$ and two corners with label $1$ lying on each side of the spine of $\omega$. In addition, these two corners are the last visited corners with label $1$ during a contour of the left side and right side of $\omega$. $F_{\infty}$ is thus delimited by the two edges joining these vertices and $v_0$, and the unique injective path in the tree joining these two vertices. The spine of $\omega$ lies in this face, except for finitely many vertices.

\bigskip

The second step takes place independently in each face of $\mathcal{M}_0$, including $F_{\infty}$. Let $F$ be a face of $\mathcal{M}_0$ and let $c_0$ be its corner with label $0$. If $F$ has finitely many vertices -- and therefore finitely many corners -- we number its corners from $0$ to $k-1$ in clockwise order along the border, starting with $c_0$. If $F$ is the infinite face, we number its corners on the right side of the spine with nonnegative integers in clockwise order, starting right after $c_0$. Similarly, we number its corners on the left side of the spine with negative integers in counterclockwise order, starting right after $c_0$. See e.g. Figure \ref{figschaeffer}. Let $\ell(i)$ denote the label of the $i$-th corner, so that $\ell(0) = 0$ and $\ell(1) = \ell(k - 1) = 1$ for a finite face whereas $\ell(1) = \ell (-1) = 1$ for $F_{\infty}$ (note that the function $\ell$ depends of the considered face).

In each face, let us define the successor function for all corners except the corners with label $0$ or $1$ by
\[ s(i) = \begin{cases}
\min \left\{ j > i : \, \ell(j) = \ell(i) - 1 \right\} & \text{if $i < 0$},\\
\min \left\{ j > i : \, \ell(j) = \ell(i) - 1 \right\} & \text{if $i>0$ and
$ \left\{ j>i : \, \ell(j) = \ell(i) -1 \right\} \neq \emptyset $},\\
\min \left\{ j \leqslant 0 : \, \ell(j) = \ell(i) - 1 \right\} & \text{if
  $i> 0$ and $ \left\{ j>i : \, \ell(j) = \ell(i) -1 \right\} = \emptyset $.}
\end{cases} \]
For a finite face, only the second case occurs, while for $F_{\infty}$ the second property of Definition \ref{spine} ensures that $\{ j \leqslant 0 : \, \ell(j) = \ell(i) - 1\}$ is finite.

\bigskip

\begin{step}
In every face, for each corner $i$ with label $\ell(i) \geqslant 2$ and such that $|s(i) - i| \neq 1$ a chord $(i, s(i))$ is added inside the face.
\end{step}

\begin{proposition}[\cite{CD}, Property 6.1]
\label{intersect}
Step 2 can be done in such a way that the various chords $(i,s(i))$ do not intersect.
\end{proposition}

\begin{remark}
The condition $|s(i) - i| \neq 1$ means that the chord $(i,s(i))$ does not already exist in $\omega$. In $F_{\infty}$, a chord $(i, s(i))$ can connect two corners that lie on different sides of the spine (see e.g. Figure \ref{figschaeffer}). This happens in the third case occurring in the definition of $s(i)$.  In that case, the corner $i$ is visited after the last occurrence of the label $\ell(i) -1$ during the contour of the right side of the spine.
\end{remark}

\bigskip

Step 2 defines a uniquely determined regular planar map $\mathcal{M}_1$ whose faces are described by the following proposition:
\begin{proposition}[\cite{CD}, Property 6.2]
\label{faces}
The faces of $\mathcal{M}_1$ are either triangular with labels $l$, $l+1$, $l+1$ or quadrangular with labels $l$, $l+1$, $l+2$, $l+1$.
\end{proposition}

\bigskip

\begin{step}
All edges of $\mathcal{M}_1$ with the same label on both ends are deleted.
\end{step}
After this last step, a unique infinite quadrangulation $\Phi(\omega)$ is obtained (see \cite{CD} for details). In addition, labels of vertices in the tree $\omega$ coincide with distances from the root of the corresponding vertices in $\Phi(\omega)$. Furthermore, the function $\Phi$ is one-to-one.

\section{Uniform infinite quadrangulations}
\label{sec:quadrangulations}

This section presents two different ways to define a uniform infinite random quadrangulation of the plane.

\subsection{Direct approach}

In \cite{Kr}, the uniform infinite quadrangulation is defined as the law of the local limit of uniformly distributed finite random quadrangulations. This limit is taken with respect to the following topology: for $Q \in \mathbf{Q}$ and $R \geqslant 0$, we denote by $B_{\mathbf{Q},R}(Q)$ the union of the faces of $Q$ that have a vertex at distance strictly smaller than $R$ from the root vertex. We may view $B_{\mathbf{Q},R}(Q)$ as a finite rooted planar map. The set $\mathbf{Q}$ is equipped with the distance
\[ d_{\mathbf{Q}} (Q_1,Q_2) = \left( 1 + \sup \left\{ R : \,
B_{\mathbf{Q},R}(Q_1) = B_{\mathbf{Q},R}(Q_2) \right\} \right)^{-1} ,\]
where the equality $B_{\mathbf{Q},R}(Q_1) = B_{\mathbf{Q},R}(Q_2)$ is in the sense of equality between two finite rooted planar maps.

Let $(\overline{\mathbf{Q}},d_{\mathbf{Q}})$ be the completion of the metric space $(\mathbf{Q},d_{\mathbf{Q}})$. Elements of $\overline{\mathbf{Q}}$ that are not finite quadrangulations are called infinite rooted quadrangulations in the sense of Krikun.

\begin{figure}[!h]
\begin{center}
\includegraphics{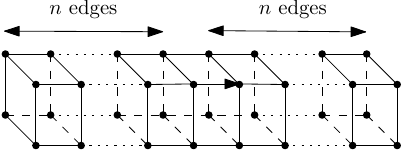}
\end{center}
\caption{A quadrangulation that converges 
in Krikun's sense to an infinite quadrangulation that is not an infinite planar map.}
\label{figexample}
\end{figure}

Note that this definition is not equivalent to Definition \ref{Qinf}. For example, the quadrangulation $Q_n$ of Figure \ref{figexample} converges as $n$ goes to infinity in $(\overline{\mathbf{Q}},d_{\mathbf{Q}})$ to an infinite quadrangulation $Q$ in Krikun's sense that is not an infinite planar map in the sense of Definition \ref{acceptable}: any proper embedding of $Q$ in $\mathbb{R}^2$ is not locally finite.

\begin{theorem}[\cite{Kr}, Theorem 1]
\label{ThKr}
For every $n \geqslant 1$ let $\nu_n$ be the uniform probability measure on
$\mathbf{Q}_n$. The sequence $(\nu_n)_{n \in \mathbb{N}}$ converges to a probability measure $\nu$ in the sense of weak convergence in the space  of all probability measures on $(\overline{\mathbf{Q}}, d_{\mathbf{Q}})$. Moreover, $\nu$ is supported on the set of infinite rooted quadrangulations (in the sense of Krikun).
\end{theorem}

\begin{remark}
One can extend the function $Q \in \mathbf{Q} \mapsto B_{\mathbf{Q},R}(Q)$ to a continuous function
$B_{\overline{\mathbf{Q}},R}$ on $\overline{\mathbf{Q}}$. $B_{\overline{\mathbf{Q}},R}(Q)$ is naturally interpreted as the union of faces of $Q$ that have a vertex at distance strictly smaller than $R$ from the root.
\end{remark}

\subsection{Indirect approach}
\label{sec:indirect}

Another possible approach to define a uniform infinite random quadrangulation is to start from a uniform infinite well-labelled tree and to consider the image of its law under Schaeffer's correspondence. This method has been developed in \cite{CD}, to which we refer for details and for proofs of what follows in this section. Let us equip $\overline{\mathbb{T}}$ with the distance
\[ d_{\mathbb{T}}(\omega,\omega') = \left( 1 + \sup \left\{ S : \,
  B_{\mathbb{T},S}(\omega) = B_{\mathbb{T},R} (\omega') \right\} \right)^{-1} ,\]
where $B_{\mathbb{T},S} (\omega)$ is the subtree of $\omega$ up to generation $S$. The metric space $(\overline{\mathbb{T}},d_{\mathbb{T}})$ is complete.

We have the following result:
\begin{theorem}[\cite{CD}, Theorem 3.1]
\label{cvarbres}
Let $\mu_n$ be the uniform probability measure on the set of all well-labelled trees with $n$ edges.
The sequence $(\mu_n)_{n \in \mathbb{N}}$ converges weakly to a probability measure $\mu$ supported on $\mathbb{T}_{\infty}$. This limit law is called the law of the uniform infinite well-labelled tree.
\end{theorem}

One of the key steps to prove this result is to show the convergence
\[ \mu_n \left( \omega \in \overline{\mathbb{T}} : B_{\overline{\mathbb{T}},S}
(\omega) = \omega^{\star} \right) \underset{n \to \infty}{\longrightarrow} \mu \left(
\omega \in \overline{\mathbb{T}} : B_{\overline{\mathbb{T}},S} (\omega) = \omega^{\star} \right) \]
for every integer $S > 0$ and every well-labelled tree $\omega^{\star}$ of height $S$. This is done by explicit computations. Let $\omega^{\star}$ be a well-labelled tree of height $S$ , and assume that $\omega^{\star}$ has exactly $k$ vertices at generation $S$, with respective labels $l_1,\ldots , l_k$. Then,
\begin{align}
\label{boulen}
\mu_n \left( \omega \in \overline{\mathbb{T}} : B_{\overline{\mathbb{T}},S}
(\omega) = \omega^{\star} \right) & = \frac{1}{D_n} \sum_{n_1 + \cdots + n_k
= n - |\omega ^{\star}|} \prod_{j=1}^k D_{n_j}^{(l_j)} ,\\
\label{boule}
\mu \left( \omega \in \overline{\mathbb{T}} : B_{\overline{\mathbb{T}},S}
(\omega) = \omega^{\star} \right) & = \frac{1}{12^{|\omega^{\star}|}} \sum_{i =
1}^k d_{l_i} \prod_{j \neq i} w_{l_j},
\end{align}
where, for every $l \geqslant 1$, $D_n^{(l)}$ is the cardinal of $\mathbb{T}_n^{(l)}$,
$D_n^{(1)} = D_n$ and
\begin{align}
\label{defw}
w_{l} & = 2 \frac{l (l +3)}{(l +1)(l +2)},\\
d_{l} & = \frac{2w_{l}}{560} (4 l^4 + 30 l^3 + 59 l^2 + 42 l
+ 4).
\label{defd}
\end{align}

\begin{proposition}[\cite{CD}, Theorem 4.3 and Theorem 5.9]
\label{etiquettes}
The measure $\mu$ is supported on $\mathscr{S}$. Furthermore,
\[\mathbb{E}_{\mu} \left[ N_l \right] =
  {\it O}(l^3) \quad \text{as $l \to \infty$}.\]
\end{proposition}

A tree with law $\mu$ has almost surely a unique spine; \cite{CD} gives a precise description of the law of the labels of this spine and of the subtrees attached to each of its vertices. For every $l> 0$, let $\rho^{(l)}$ be the measure on $\mathbf{T}^{(l)}$ defined by $\rho^{(l)}(\omega) = 12^{-|\omega|}$ for every $\omega \in \mathbf{T}^{(l)}$. Then $\frac{1}{2} \rho^{(l)}$ is the law of the Galton-Watson tree with geometric offspring distribution with parameter $\frac{1}{2}$ and with random labels generated according to the following rules. The root has label $l$ and the label of every other vertex is chosen uniformly in $\{m-1,m,m+1\}$ where $m$ is the label of its parent. Furthermore, these choices are made independently for every vertex. Proposition 2.4 of \cite{CD} proves that $\rho^{(l)}(\mathbb{T}^{(l)}) = w_l$, therefore the measure $\widehat{\rho}^{(l)}$ defined on $\mathbb{T}^{(l)}$ by $\widehat{\rho}^{(l)}(\omega) = w_l^{-1} \rho^{(l)}(\omega) =w_l^{-1} 12^{-|\omega|}$ for every $\omega \in \mathbb{T}^{(l)}$ is a probability measure. The following result will be useful for our purposes.
\begin{theorem}[\cite{CD}, Theorem 4.4]
\label{descmu}
Let $\omega$ be a random tree distributed according to $\mu$ and let $u_0, u_1, u_2, \ldots $ be the sequence of the vertices of its spine listed in genealogical order. For every $n \geqslant 0$, let $Y_n$ be the label of $u_n$.
\begin{enumerate}
\item The process $(Y_n)_{n \geqslant 0}$ is a Markov chain taking values in $\mathbb{N}$ with transition kernel $\Pi$ defined by:
\begin{align*}
\Pi (l,l-1) & = \frac{(w_l)^2}{12 d_l} d_{l-1} := q_l& \text{if $l
\geqslant 2$,}\\
\Pi (l,l) & = \frac{(w_l)^2}{12} := r_l& \text{if $l \geqslant 1$,}\\
\Pi (l,l+1) & = \frac{(w_l)^2}{12 d_l} d_{l+1} := p_l& \text{if $l \geqslant 1$.}
\end{align*}
\item Conditionally given $(Y_n)_{n \geqslant 0} = (y_n)_{n
\geqslant 0}$, the sequence $(L_n)_{n \geqslant 0}$ of subtrees of $\omega$ attached to the left side of the spine and the sequence $(R_n)_{n \geqslant 0}$ of subtrees attached to the right side of the spine form two independent sequences of independent labelled trees distributed according to the measures $\widehat{\rho}^{(y_n)}$.
\end{enumerate}
\end{theorem}

We can now map the law of the uniform infinite random tree on the set of quadrangulations using Schaeffer's correspondence. Let us equip $\Phi(\mathscr{S})$ with the distance $d_{\Phi}$ so that $\Phi$ is an isometry from $\mathscr{S}$ onto $\Phi(\mathscr{S})$. We denote by $\mu_{\Phi,n}$ and $\mu_{\Phi}$ the respective image measures of $\mu_n$ and $\mu$ under $\Phi$. The measure $\mu_{\Phi}$ is well defined because $\mu$ is supported on $\mathscr{S}$.

Since $\Phi$ is a bijection between $\mathbb{T}_n$ and $\mathbf{Q}_n$, $\mu_{\Phi,n} = \nu_n$ is the uniform probability measure on the set of quadrangulations with $n$ faces. As a direct consequence of Theorem \ref{cvarbres}, the sequence $(\mu_{\Phi,n})_{n \in \mathbb{N}}$ converges weakly to $\mu_{\Phi}$ in the space of all probability measures on $(\Phi(\mathscr{S}),d_{\Phi})$. Thus, in some sense, $\mu_{\Phi}$ can also be viewed as a uniform probability measure on the space of infinite quadrangulations.

\begin{remark}
The topology induced by $d_{\Phi}$ on the set $\Phi(\mathscr{S})$ is rather different than the one that would be induced by $d_{\mathbf{Q}}$. Indeed it may happen that two trees $\omega$ and $\omega'$ are close for the metric $d_{\mathbb{T}}$, but the quadrangulations $\Phi(\omega)$ and $\Phi(\omega')$ are very different for $d_{\mathbf{Q}}$. For example, the linear tree $\omega_n$ with $2n-1$ vertices and with labels given by the sequence $1, 2,\ldots,n-1, n , n-1, \ldots ,2,1$ converges as $n$ goes to infinity to the infinite linear tree $\omega$ with labels given by the sequence $1,2, \ldots$.
As a consequence, the quadrangulation $\Phi(\omega_n)$ converges to the infinite quadrangulation $\Phi(\omega)$ in $(\Phi(\mathscr{S}),d_{\Phi})$ as $n$ goes to infinity. On the other hand, for every $n \geqslant 1$, the quadrangulation $\Phi(\omega_n)$ has two vertices at distance $1$ from its root whereas $\Phi(\omega)$ has only one vertex at distance $1$ from its root and therefore the sequence $(\Phi(\omega_n))_{n \in \mathbb{N}}$ does not converge to $\Phi(\omega)$ in $(\overline{\mathbf{Q}},d_{\mathbf{Q}})$.

It is then a natural question to ask whether the two notions of uniform infinite quadrangulation that we have introduced coincide.
\end{remark}

\section{Equality of the two uniform infinite quadrangulations}
\label{sec:egalite}

In this section, we will show that the two definitions of the uniform infinite quadrangulation coincide. The first problem comes from the fact that we have two different notions of infinite quadrangulations: elements of $\Phi(\mathscr{S})$, which are regular planar maps on one hand, and elements of the completion $\overline{\mathbf{Q}}$ of $\mathbf{Q}$ on the other hand. This problem can be solved by identifying $\Phi(\mathscr{S})$ with a subset of $\overline{\mathbf{Q}}$, allowing us to consider $\mu_{\Phi}$ as a measure on $\overline{\mathbf{Q}}$ supported on $\Phi(\mathscr{S})$.

More precisely, let $R>0$ and $\omega \in \mathscr{S}$. Define $B_R(\Phi(\omega))$ as the union of all faces of $\Phi(\omega)$ that have a vertex at distance strictly smaller than $R$ from the root. Since the tree $\omega$ has only finitely many vertices with label smaller than $R+1$, there are finitely many such faces and $B_R(\Phi(\omega))$ is a finite map. Therefore $\mathbb{S}^2 \setminus B_R(\Phi(\omega))$ has finitely many connected components; and the boundaries of these components are finite length cycles of $\Phi(\omega)$.

Let $\gamma$ be such a cycle. Each edge of $\gamma$ is adjacent to two faces of $\Phi(\omega)$. One has a vertex at distance strictly smaller than $R$ from the root, and the other one has only vertices at distance at least $R$ from the root. The quadrangulation being bipartite, each edge of $\gamma$ connects a vertex at distance $R$ from the root with a vertex at distance $R+1$ from the root. Therefore, by adding to $B_R(\Phi(\omega))$ an extra vertex in the connected component of $\mathbb{S}^2 \setminus B_R(\Phi(\omega))$ bounded by $\gamma$ and an edge between this vertex and each vertex of $\gamma$ at distance $R+1$ from the root, and repeating this operation for every connected component of $\mathbb{S}^2 \setminus B_R(\Phi(\omega))$, we obtain a finite quadrangulation. The sequence of finite quadrangulations obtained in this way for every $R>0$ converges to $\Phi(\omega)$ as $R$ goes to infinity, in the sense of Krikun, showing that for every tree $\omega \in \mathscr{S}$, $\Phi(\omega)$ can be identified with an element of $\overline{\mathbf{Q}}$.

To be able to consider $\mu_{\Phi}$ as a measure on $\overline{\mathbf{Q}}$, we now need to verify that the mapping $\Phi : \mathscr{S} \to \overline{\mathbf{Q}}$ is measurable with respect to the Borel $\sigma$-field of $\left( \overline{\mathbf{Q}}, d_{\mathbf{Q}} \right)$. The following lemma is proved in Section \ref{propbij}:
\begin{lemma}
\label{measurable}
Fix $R > 0$ and $\omega_0 \in \mathscr{S}$. The set $ A = \{ \omega \in \mathscr{S} : \,  B_{\overline{\mathbf{Q}},R}(\Phi(\omega)) = B_{\overline{\mathbf{Q}},R}(\Phi(\omega_0)) \}$ is measurable with respect to the Borel $\sigma$-field of $\left( \mathscr{S}, d_{\overline{\mathbb{T}}} \right)$.
\end{lemma}

Fix $Q^{\star} \in \overline{\mathbf{Q}}$. Lemma \ref{measurable} implies that
\begin{equation*}
 \Phi^{-1} \left(
\left\{ Q \in \overline{\mathbf{Q}} : \, d_{\mathbf{Q}} (Q,Q^{\star}) \leqslant \frac{1}{R+1} \right\}                 \right)
= \Phi^{-1} \left(
\left\{ \vphantom{\frac{1}{R}} Q \in \overline{\mathbf{Q}} : \, B_{\overline{\mathbf{Q}},R}(Q) =
                     B_{\overline{\mathbf{Q}},R}(Q^{\star}) \right\}
\right)
\end{equation*}
is measurable with respect to the Borel $\sigma$-field of $\left( \mathscr{S} , d_{\overline{\mathbb{T}}} \right)$, proving that $\Phi : \mathscr{S} \to \overline{\mathbf{Q}}$ is measurable. Therefore, we may and will see $\mu_{\Phi}$ as a probability measure on $\left( \overline{\mathbf{Q}} , d_{\mathbf{Q}} \right)$.

\bigskip

We are now ready to state our main result:
\begin{theorem}
\label{equality}
The sequence $\left( \mu_{\Phi,n} \right)_{n \in \mathbb{N}}$ converges weakly to $\mu_{\Phi}$ in the space of all probability measures on $\left(\overline{\mathbf{Q}}, d_{\mathbf{Q}} \right)$. Therefore $\mu_{\phi}$ viewed as a probability measure on $\left(\overline{\mathbf{Q}}, d_{\mathbf{Q}} \right)$ coincides with $\nu$.
\end{theorem}

Since $\mu_{\Phi,n} = \nu_n$ and $\nu$ is defined as the limit of the sequence $(\nu_n)$ in the space of all probability measures on $\left(\overline{\mathbf{Q}}, d_{\mathbf{Q}} \right)$, the second assertion is a direct consequence of the first one.

To establish the first assertion, we have to show that for every $Q^{\star} \in \overline{\mathbf{Q}}$ and $R > 0$ one has
\begin{equation*}
\mu_{\Phi,n} \Big( Q \in \overline{\mathbf{Q}} :
B_{\overline{\mathbf{Q}},R}(Q) = B_{\overline{\mathbf{Q}},R}(Q^{\star})\Big)
\underset{n \to \infty}{\longrightarrow}
\mu_{\Phi} \Big( Q \in \overline{\mathbf{Q}} :
B_{\overline{\mathbf{Q}},R}(Q) = B_{\overline{\mathbf{Q}},R}(Q^{\star})\Big).
\end{equation*}
The remaining part of this work is devoted to the proof of this convergence.

\subsection{A property of Schaeffer's correspondence}
\label{propbij}

For every integers $S>0$ and $R>0$ we let
\begin{equation}
\label{omegas}
\Omega_S (R) = \{ \omega \in \overline{\mathbb{T}} : \,
\omega \text{ has a label} \leqslant R+1 \text{ strictly above generation } S \}.
\end{equation}
In the first two statements of this section, $S$ and $R$ are two fixed positive integers.

\begin{proposition}
\label{facegen}
Let $\omega$ be a tree of $\mathscr{S}$ which does not belong to $\Omega_S(R)$ (i.e. $\omega$ is such that the label of every vertex at a generation strictly greater than $S$ is at least $R+2$). Then $B_{\overline{\mathbf{Q}},R}(\Phi(\omega)) = B_{\overline{\mathbf{Q}},R} (\Phi(B_{\overline{\mathbb{T}},S}(\omega)))$.
\end{proposition}
\begin{proof}
The proof follows step by step the construction of $\Phi(\omega)$ in Section \ref{subsec:schaeffer}. Fix an embedding of $\omega$ as an infinite planar map.

In the first step, an infinite planar map $\mathcal{M}_0(\omega)$ is obtained from $\omega$ by adding an extra vertex $v_0$ with label $0$ and edges between $v_0$ and corners with label $1$. Similarly, we can construct a planar map $\mathcal{M}_0(B_{\overline{\mathbb{T}},S}(\omega))$. The extra edges in these two maps are uniquely determined by corners with label $1$, and these corners are determined by $B_{\overline{\mathbb{T}},S}(\omega)$ (no vertex at a generation greater than $S$ has a label less than $2$). We consider the unique ``infinite face'' of $\mathcal{M}_0$ as an extra face. The maps $\mathcal{M}_0(\omega)$ and $\mathcal{M}_0(B_{\overline{\mathbb{T}},S}(\omega))$ then have the same number of faces, say $p$, which in addition have the same boundaries, composed by the two edges joining $v_0$ to the corners with label $1$ and, in the case when these corners with label $1$ belong to different vertices, the unique injective path in the tree between these two vertices. Let $F_1(\omega),\ldots, F_p(\omega)$ and $F_1(B_{\overline{\mathbb{T}},S}(\omega)),\ldots , F_p(B_{\overline{\mathbb{T}},S}(\omega))$ denote the faces of $\mathcal{M}_0(\omega)$ and $\mathcal{M}_0(B_{\overline{\mathbb{T}},S}(\omega))$ respectively, listed in such a way that, for every $i$, the faces $F_i(\omega)$ and $F_i(B_{\overline{\mathbb{T}},S}(\omega))$ have the same boundary.

In the second step, edges $(c,s(c))$ are added inside each face for every corner $c$, finally giving two regular planar maps $\mathcal{M}_1(\omega)$ and $\mathcal{M}_1(B_{\overline{\mathbb{T}},S}(\omega))$. Let us consider a face $F_i(\omega)$ of $\mathcal{M}_0(\omega)$ and the corresponding face $F_i(B_{\overline{\mathbb{T}},S}(\omega))$. The corners of these faces are numbered $(c_{i,j})_{j \in J_i}$ for
$F_i(\omega)$ and $(c'_{i,j})_{j \in J'_i}$ for $F_i(B_{\overline{\mathbb{T}},S}(\omega))$, the numbering being in clockwise order for a finite face and counterclockwise order for corners on the left side if the spine of the tree, clockwise order for corners on the right side of the spine of the tree, in the case of the infinite face.

For every $i \in \{1, \ldots ,p \}$, let $v_{i,1}, \ldots , v_{i,k_i}$ be the vertices of $F_i(\omega)$ at generation $S$ that have at least one child. These vertices are also vertices of $F_i(B_{\overline{\mathbb{T}},S}(\omega))$ at generation $S$ and their labels are greater than $R+1$. For every $j \leqslant k_i$, let $e_{i,j}$ be the last corner  before $v_{i,j}$ in $F_i(\omega)$ and with label $R+1$. This corner is the same in $F_i(\omega)$ and $F_i(B_{\overline{\mathbb{T}},S}(\omega))$. The same edge $(e_{i,j},s(e_{i,j}))$ joining $e_{i,j}$ to the first corner following $e_{i,j}$ with label $R$ is thus added to $F_i(\omega)$ and $F_i(B_{\overline{\mathbb{T}},S}(\omega))$ (note that this corner is also the first corner with label $R$ following every corner of $v_{i,j}$, see Figure \ref{faceisole}).

\begin{figure}[!t]
\begin{center}
\includegraphics[width=0.6\textwidth]{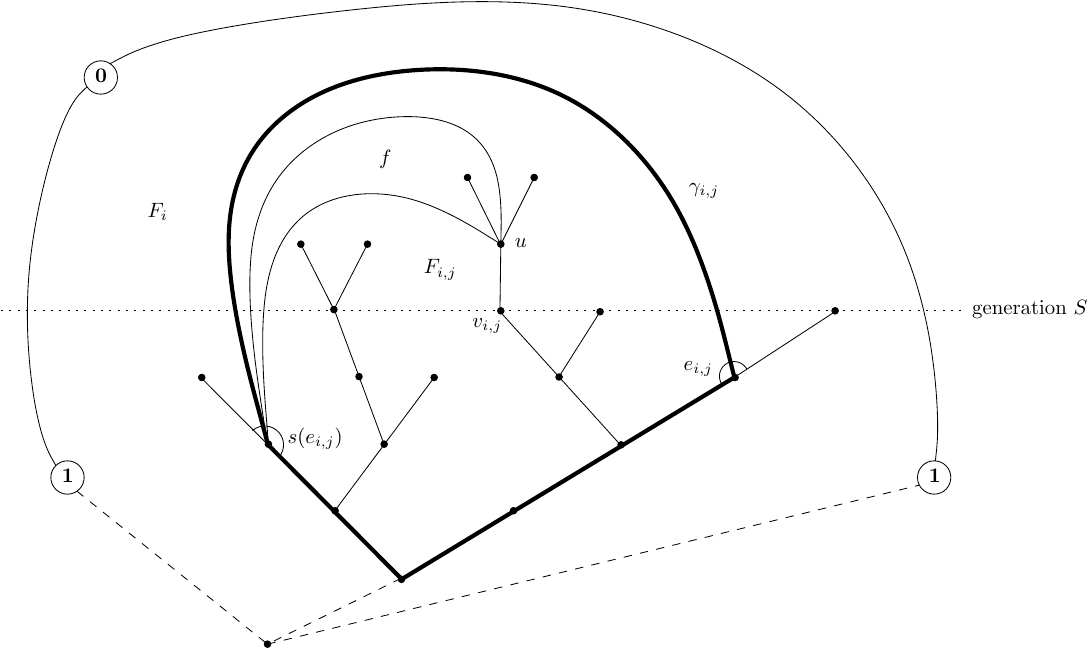}
\end{center}
\caption{A face $F_i$ and a cycle $\gamma_{i,j}$ associated with a vertex $v_{i,j}$ at generation $S$.}
\label{faceisole}
\end{figure}

Therefore for every $j \in \{ 1, \ldots k_i \}$ the same cycle $\gamma_{i,j}$ composed by the edge $(e_{i,j},s(e_{i,j}))$ and the genealogical path between $e_{i,j}$ and $s(e_{i,j})$ appears in both $F_i(\omega)$ and $F_i(B_{\overline{\mathbb{T}},S}(\omega))$ (see Figure \ref{faceisole}). For $j \neq j'$ the (strict) interiors of the cycles $\gamma_{i,j}$ and $\gamma_{i,j'}$ are either disjoint, or one of them is contained in the other one. Here we define the interior of a cycle as the connected component of the complement of this cycle which does not contain $v_0$.

Let us now show that if a face $f$ of $\Phi(\omega)$ intersects the interior of a cycle $\gamma_{i,j}$, the labels of vertices of $f$ are greater than or equal to $R$.  We first deal with the case when $f$ has a vertex $u$ that belongs to the interior of the cycle $\gamma_{i,j}$. If the label of $u$ is greater than or equal to $R+2$ the conclusion is obvious. If not, the label of $u$ is $R+1$ and for $f$ to have a vertex with label $R-1$, $u$ must be connected to vertices with label $R$ by two edges: the only possible choice for a vertex with label $R$ is $s(e_{i,j})$ and the last vertex of $f$ would then belong to the domain bounded by the union of the two edges connecting $u$ to $s(e_{i,j})$ (see Figure \ref{faceisole}) so that its label could not be $R-1$. The case when no vertex of $f$ belongs to the interior of the cycle $\gamma_{i,j}$ is treated in a similar manner.

\bigskip

The previous discussion shows that faces of $\Phi(\omega)$, respectively of $\Phi(B_{\overline{\mathbb{T}},S}(\omega))$, that intersect the interior of a cycle $\gamma_{i,j}$, are not taken into account in the definition of $B_{\overline{\mathbf{Q}},R}(\Phi(\omega))$, respectively of $B_{\overline{\mathbf{Q}},R}(\Phi(B_{\overline{\mathbb{T}},S}(\omega)))$.

Let us denote by $\Phi_1(\omega)$ the planar map obtained after the second step of the construction of $\Phi(\omega)$ (this map is denoted by $\mathcal{M}_1$ in Section \ref{subsec:schaeffer}). Let us consider the map $\widetilde{\Phi}_1(\omega)$ obtained by removing every edge and vertex of $\Phi_1(\omega)$ lying in the interior of a cycle $\gamma_{i,j}$. By construction, every vertex of $\omega$ with generation strictly greater than $S$ belongs to the interior of a cycle $\gamma_{i,j}$. It follows that
\begin{equation}
\label{phi1}
\widetilde{\Phi}_1(\omega) = \widetilde{\Phi}_1 \left( B_{\overline{\mathbb{T}},S} (\omega) \right).
\end{equation}

Finally, let $\Phi_2(\omega)$ denote the map obtained by removing every edge of $\widetilde{\Phi}_1(\omega)$ connecting two vertices of $\widetilde{\Phi}_1(\omega)$ with the same label less than or equal to $R$. Every face of $\Phi(\omega)$ that is taken into account in the ball $B_{\overline{\mathbf{Q}},R}(\Phi(\omega))$ is also a quadrangular face of $\Phi_2(\omega)$. Conversely, every quadrangular face of $\Phi_2(\omega)$ having a vertex with label strictly smaller than $R$ is also a face of $B_{\overline{\mathbf{Q}},R}(\Phi(\omega))$. In other words, the ball $B_{\overline{\mathbf{Q}},R}(\Phi(\omega))$ is the union of the quadrangular faces of $\Phi_2(\omega)$ having a vertex with label strictly smaller than $R$. From \eqref{phi1} we have
\[
\Phi_2(\omega) = \Phi_2 \left(B_{\overline{\mathbb{T}},S} (\omega) \right),
\]
and the previous observations allow us to conclude that
\[
B_{\overline{\mathbf{Q}},R} \left( \Phi (\omega) \right) =
B_{\overline{\mathbf{Q}},R} \left( \Phi \left(
B_{\overline{\mathbb{T}},S} (\omega) \right) \right)
\]
which completes the proof.
\end{proof}

\begin{corollary}
\label{pbcle}
Let $\omega_0 \in \mathscr{S}$. There exists a countable collection $\left(\omega^{S,R}_i\right)_{i \in I}$ of trees in $\mathscr{S} \cap \Omega_S(R)^c$ verifying for all $i \in I$
\[
B_{\overline{\mathbf{Q}},R} \left( \Phi (\omega_i^{S,R}) \right) = B_{\overline{\mathbf{Q}},R} \left( \Phi (\omega_0) \right)
\]
and such that for every $\omega \in \mathscr{S} \cap \Omega_S(R)^c$ the following assertions are equivalent:
\begin{enumerate}
 \item $B_{\overline{\mathbf{Q}},R} \left( \Phi (\omega) \right) = B_{\overline{\mathbf{Q}},R} \left( \Phi (\omega_0) \right);$
 \item there exists $i \in I$ such that $B_{\overline{\mathbb{T}},S} (\omega) = B_{\overline{\mathbb{T}},S} (\omega_i^{S,R})$.
\end{enumerate}
\end{corollary}
\begin{proof}
The collection $\left(\omega_i^{S,R}\right)_{i \in I}$ that consists of all finite trees $\omega'$ having at most $S$ generations and such that $B_{\overline{\mathbf{Q}},R} \left( \Phi (\omega') \right) = B_{\overline{\mathbf{Q}},R} \left( \Phi (\omega_0) \right)$ is countable and has the desired properties. Indeed, if $\omega \in \mathscr{S} \cap \Omega_S(R)^c$ and if there exists $i \in I$ such that $B_{\overline{\mathbb{T}},S} (\omega) = B_{\overline{\mathbb{T}},S} (\omega_i^{S,R})$, Proposition \ref{facegen} ensures that $B_{\overline{\mathbf{Q}},R} \left( \Phi (\omega) \right) = B_{\overline{\mathbf{Q}},R} \left( \Phi (\omega_i^{S,R}) \right) = B_{\overline{\mathbf{Q}},R} \left( \Phi (\omega_0) \right)$. Conversely, if $B_{\overline{\mathbf{Q}},R} \left( \Phi (\omega) \right) = B_{\overline{\mathbf{Q}},R} \left( \Phi (\omega_0) \right)$ and $\omega \in \mathscr{S} \cap \Omega_S(R)^c$, then $\omega' = B_{\overline{\mathbb{T}},S} (\omega)$ verifies $B_{\overline{\mathbf{Q}},R} \left( \Phi (\omega) \right) = B_{\overline{\mathbf{Q}},R} \left( \Phi (\omega') \right)$ by Proposition \ref{facegen} and $\omega'$ belongs to the collection $\left(\omega_i^{S,R}\right)_{i \in I}$.
\end{proof}

We conclude this section with the proof of Lemma \ref{measurable}.
\begin{proof}[Proof of Lemma \ref{measurable}]
Fix $R>0$. For every $S>0$ the set $\Omega_S(R)$ is open and closed in $\overline{\mathbb{T}}$. In addition one has
\begin{equation*}
\begin{split}
A & = \bigcup_{S > 0} \left( \left\{ \omega
\in \mathscr{S} : \, B_{\overline{\mathbf{Q}},R} (\Phi(\omega)) =
B_{\overline{\mathbf{Q}},R} (\Phi(\omega_0)) \right\} \cap \Omega_S(R)^c \right)\\
& = \bigcup_{S > 0} \bigcup_{i \in I_{S,R}} \left( \left\{ \omega \in
  \mathscr{S} : \, B_{\overline{\mathbb{T}},S} (\omega) = B_{\overline{\mathbb{T}},S}
(\omega_i^{S,R}) \right\} \cap \Omega_S(R)^c \right)
\end{split}
\end{equation*}
where $\left(\omega_i^{S,R}\right)_{i \in I_{S,R}}$ is the collection given by Corollary \ref{pbcle}. This shows that the set $A$ is measurable.
\end{proof}

\subsection{Asymptotic behavior of labels on the spine}
\label{spinelabels}

Recall that the sequence $(Y_k)_{k \geqslant 0}$ of the successive labels of vertices of the spine is a Markov chain with transition matrix $\Pi$ given by Theorem \ref{descmu}. In this section we study the asymptotic behavior of this Markov chain.
\begin{lemma}
\label{labspine}
The Markov chain $(Y_k)_{k \geqslant 0}$ is transient. In addition, for every $\varepsilon > 0$ there exists $\alpha > 0$ such that for $k$ large enough one has
\[
\mathbb{P}\left[ Y_j \geqslant \alpha k, \, \forall j \geqslant 0 \, \middle|
  \, Y_0   =k \right] \geqslant 1 - \varepsilon.
\]
\end{lemma}

\begin{proof}
The Taylor expansion $\frac{q_k}{p_k} = 1 - \frac{8}{k} +
{\it O} \left( \frac{1}{k^2} \right)$ (\cite{CD}, Lemma 5.5) implies that there exists $C >0$ such that
\[
\prod_{i=2}^{k} \frac{q_i}{p_i} \underset{k \to \infty}{\sim} C k^{-8}.
\]
A standard argument for birth and death processes then ensures that $Y$ is transient.  Furthermore, for every $k > j \geqslant 1$,
\[
\mathbb{P}_k \left[T_j = \infty \right] =
  \frac{\sum_{i=j}^{k-1}\frac{q_i}{p_i} \frac{q_{i-1}}{p_{i-1}} \cdots \frac{q_{j+1}}{p_{j+1}}}
       {\sum_{i=j}^{\infty}\frac{q_i}{p_i} \frac{q_{i-1}}{p_{i-1}} \cdots \frac{q_{j+1}}{p_{j+1}}}
\]
where $T_j$ is the hitting time of $j$. Therefore one has, for $\alpha <
1$,
\begin{align*}
\mathbb{P}_k \left[T_{[\alpha k]} = \infty \right] & =
  \frac{\frac{1}{k}\sum_{i=[\alpha k]}^{k-1}\frac{q_i}{p_i} \frac{q_{i-1}}{p_{i-1}} \cdots
    \frac{q_{[\alpha k] +1}}{p_{[\alpha k] +1}}}
       {\frac{1}{k}\sum_{i=[\alpha k]}^{\infty}\frac{q_i}{p_i} \frac{q_{i-1}}{p_{i-1}} \cdots
         \frac{q_{[\alpha k]+1}}{p_{[\alpha k]+1}}}\\
& \underset{k \to \infty}{\longrightarrow} 
\frac{\int_{\alpha}^1 \left( \frac{\alpha}{t} \right)^8 dt}
     {\int_{\alpha}^{\infty} \left( \frac{\alpha}{t} \right)^8 dt} =
     1 - \alpha^7.
\end{align*}
The desired result follows.
\end{proof}

\begin{proposition}
\label{limlabspine}
Let $Z$ be a nine-dimensional Bessel process started at $0$. Then
\[
 \left( \frac{1}{\sqrt{n}} Y_{[nt]} \right)_{t \geqslant 0} \underset{n \to
   \infty}{\longrightarrow} \left( Z_{\frac{2}{3} t} \right)_{t \geqslant 0}
\]
in the sense of convergence in distribution in the space $D(\mathbb{R}_+,\mathbb{R}_+)$.
\end{proposition}
\begin{proof}
The convergence in the proposition is a direct consequence of a more general result by Lamperti \cite{Lam} which we now recall. Let $(X_n)_{n \geqslant 0}$ be a time-homogeneous Markov chain on $\mathbb{R}_+$ verifying:
\begin{enumerate}
\item for every $K > 0$ one has uniformly in $x \in \mathbb{R}_+$
\[\lim_{n \rightarrow \infty} \frac{1}{n} \sum_{i = 0}^{n-1}
\mathbb{P} \left( X_i \leqslant K \, \middle| \, X_0 = x \right) = 0;\]
\item for every $k \in \mathbb{N}$ the following moments exist and are bounded as functions of $x \in \mathbb{R}_+$
\[ m_k (x) = \mathbb{E} \left[ (X_{n+1} - X_n)^k \middle| X_n = x \right];\]
\item there exist $\beta > 0$ and $\alpha > - \beta /2$ such that
\begin{align*}
\lim_{x \rightarrow \infty} m_2(x) & = \beta,\\
\lim_{x \rightarrow \infty} x \, m_1(x) & = \alpha.
\end{align*}
\end{enumerate}
Let us define the process $\left(x_t^{(n)} \right)_{t \in \mathbb{R}_+}$ by $x_t^{(n)} = n^{-1/2} X_i$
if $t = \frac{i}{n}$, $i = 0,1,2, \ldots$, and linear interpolation on intervals of the form $\left[ \frac{i-1}{n}, \frac{i}{n} \right]$. Lamperti's theorem states that $\left(x_t^{(n)}\right)_{t \in \mathbb{R}_+}$ converges in distribution to the diffusion process $(x_t)_{t \in \mathbb{R}_+}$ with generator
\[ L = \frac{\alpha}{x} \frac{\mathrm{d}}{\mathrm{d}x} + \frac{\beta}{2}
\frac{\mathrm{d}^2}{\mathrm{d}x^2}.\]

In our case, we consider the Markov chain $\widetilde{Y}$ whose transition matrix is given by $\widetilde{\Pi}
(x,y) = \Pi ([x],[y])$ if $y = x+1,x-1$ or $x$. Assertion 1. easily follows from Lemma
\ref{labspine} and assertion 2. is trivial.
In addition one has $p_n = \frac{1}{3} +
\frac{4}{3n} + {\it O} (n^{-2})$ and $q_n = \frac{1}{3} -
\frac{4}{3n} + {\it O} (n^{-2})$ (\cite{CD}, Lemma 5.5) giving:
\begin{align*}
\lim_{x \rightarrow \infty} m_2(x) & = \frac{2}{3},\\
\lim_{x \rightarrow \infty} x \, m_1(x) & = \frac{8}{3}.
\end{align*}
and therefore assertion 3. holds with $\alpha = 8/3$ and $\beta = 2/3$.

The rescaled chain $Y$ thus converges in law to the diffusion process with generator
\[L = \frac{2}{3} \left(\frac{4}{x} \frac{\mathrm{d}}{\mathrm{d}x}
+ \frac{1}{2} \frac{\mathrm{d}^2}{\mathrm{d}x^2} \right),\]
which was the desired result.
\end{proof}

\subsection{Asymptotic properties of small labels}
\label{smalllabels}

Thanks to Corollary \ref{pbcle}, the proof of Theorem \ref{equality} will reduce to showing that the $\mu_n$-measure of certain balls in the space of trees converges to the corresponding $\mu$-measure. Still we need to show that the error made by disregarding trees that belong to $\Omega_S(R)$ is small when $S$ is large.
In this section we fix $R,\varepsilon >0$ and we write $\Omega_S = \Omega_S(R)$ to simplify notation.

\begin{lemma}
\label{seuil}
There exists an integer $S^{\star} > 0$ such that $\mu(\Omega_S) \leqslant \varepsilon$ for every $S > S^{\star}$.
\end{lemma}
\begin{proof}
Let  $\Omega = \bigcap_{S=1}^{\infty} \Omega_S$. If $\omega \in \Omega$ then $\omega$ has infinitely many vertices with label in $\{1, \ldots, R+1 \}$ and there exists $l \in\{1, \ldots, R+1 \}$ such that $N_l(\omega) = \infty$. Since $\mu$ is supported on $\mathscr{S}$, one has $\mu ( \Omega) = 0$ and therefore $\mu (\Omega_S) \rightarrow 0$ as $S \rightarrow \infty$.
\end{proof}

The main ingredient of the proof of Theorem \ref{equality} is Proposition \ref{nseuil}, which gives an analog of Lemma \ref{seuil} when $\mu$ is replaced by $\mu_n$, with \emph{uniformity} in $n$. To establish this estimate, we will need an upper bound for the probability that there exists a vertex at generation $S$ with label smaller than $S^{\alpha}$, where $\alpha < 1/2$ is fixed (see Lemma \ref{estimlabels} below). Let us first give an easy preliminary lemma:
\begin{lemma}
\label{size}
Fix $S>0$. There exist positive integers $N_1(S)$ and $K_{\varepsilon}(S)$ such that for every $n > N_1(S)$:
\[ \mu_n \left( \omega : \, \left|B_{\overline{\mathbb{T}},S}(\omega)\right|
  > K_{\varepsilon}(S)  \right) < \varepsilon .\]
\end{lemma}
\begin{proof}
This result is a direct consequence of the convergence of the measures $\mu_n$ to $\mu$. Indeed, as $\left|B_{\overline{\mathbb{T}},S}(\omega)\right|$ is finite for every tree $\omega$, one can chose $K_{\varepsilon}(S)$ large enough such that
\[ \mu \left( \omega : \, \left|B_{\overline{\mathbb{T}},S}(\omega)\right| > K_{\varepsilon}(S)  \right) < \varepsilon.\]
The convergence of $\mu_n$ to $\mu$ then gives $N_1(S)$ such that the inequality of the lemma is true for $n > N_1(S)$.
\end{proof}

There exists a finite number of well-labelled trees with height exactly $S$ and having at most $K_{\varepsilon}$ edges. Let us denote this number by $M_{\varepsilon}(S)$.

\bigskip

For every $S > 0$ and $\alpha \in \left[0, \frac{1}{2} \right[$ we let
\[A_{\alpha}(S) = \left\{ \omega \in \overline{\mathbb{T}} : \,
\omega \text{ has a vertex at generation $S$ with label } \leqslant S^{\alpha}
\right\}.\]

\begin{lemma}
\label{estimlabels}
Fix $\alpha < \frac{1}{2}$. For every sufficiently large integer $S$, there exists $N_2(S)$ such that, for every $n > N_2(S)$, one has
\[
\mu_n \left( A_{\alpha}(S) \right) < \varepsilon.
\]
\end{lemma}
\begin{proof}
We first observe that it is enough to prove the bound $\mu \left( A_{\alpha}(S) \right) < \varepsilon$ when $S$ is large. Indeed the set $A_{\alpha}(S)$ is closed in $\overline{\mathbb{T}}$, and thus we have $\limsup_n \mu_n(A_{\alpha}(S)) \leqslant \mu(A_{\alpha}(S))$.

Recall the notation $\rho^{(l)}$ and $\widehat{\rho}^{(l)}$ introduced in Section \ref{sec:indirect}.
For $H > 0$ and $l>0$ one has
\begin{equation*}
\widehat{\rho}^{(l)} \left( \vphantom{\rho^l} h(\omega) > H \right) =
\frac{1}{w_l} \sum_{\substack{\omega \in \mathbb{T}^{(l)}\\ h(\omega) > H}} 12^{-|\omega|}
\leqslant
\frac{1}{w_l} \sum_{\substack{\omega \in \mathbf{T}^{(l)}\\ h(\omega) > H}} 12^{-|\omega|}
= \frac{1}{w_l} \rho^{(l)} \left( \vphantom{\rho^l} h(\omega) > H \right).
\end{equation*}
Therefore
\[
\widehat{\rho}^{(l)} \left( \vphantom{\rho^l} h(\omega) > H \right)
\leqslant \frac{2}{w_l} \mathbb{P}_{GW(1/2)} \left[\vphantom{\rho^l} h(\omega) > H \right],
\]
where $\mathbb{P}_{GW(1/2)}$ is the law of a Galton-Watson tree whose offspring distribution is geometric with parameter $1/2$. Theorem 1 (page 19) of \cite{AN} gives
\[
\lim_{H \rightarrow \infty} H \, \mathbb{P}_{GW(1/2)} \left[
  \vphantom{\rho^k} h(\omega) > H \right] = 1.
\]
From the explicit formula for $w_l$ we have $\frac{2}{w_l} \leqslant \frac{3}{2}$ for every $l \geqslant 0$. Hence there exists $H_1 > 0$ such that for $H>H_1$
\[
\widehat{\rho}^{(l)} \left(\vphantom{\rho^l}  h(\omega) > H \right) \leqslant \frac{2}{H}.
\]
Fix $\eta \in \left]0, \frac{1}{2}\right[$. Recall that $g_S(\omega)$ is the set of vertices of $\omega$ at generation $S$ and that, for every integer $k$, $L_k$ and $R_k$ are the subtrees of $\omega$ attached respectively to the left side and to the right side of the $k$-th vertex of the spine of $\omega$. For $S > (1- \eta)^{-1} H_1$, Theorem \ref{descmu} and the previous bound give
\[
\mu \left[ g_S(\omega) \cap \bigcup_{0 \leqslant k \leqslant
    [\eta S] - 1} (L_k \cup R_k) \neq \emptyset \right] \leqslant 2
\sum_{k=1} ^{[\eta S] - 1} \frac{2}{S-k} \leqslant 4 \frac{\eta}{1
  - \eta} \leqslant 8 \eta,
\]
and therefore
\[
\mu \left( A_{\alpha}(S) \right) \leqslant 8 \eta +
\mu \left( \exists s \in g_S(\omega) \cap \bigcup_{k = [\eta S]}^{S}
\left( L_k \cup R_k \right) \, : \, \ell(s) \leqslant S^{\alpha} \right).
\]
Applying the Markov property at time $[\eta S]$ to the Markov chain $Y$ and then using Proposition \ref{limlabspine} and Lemma \ref{labspine} we find $\delta > 0$ and $S_1$ such that for $S > S_1$ one has
\[
\mu \left( Y_k \geqslant [\delta \sqrt{S}] , \, \forall k \geqslant
  [\eta S] \right) \geqslant 1 - \eta.
\]
We now have
\begin{equation}
\label{Aalphaintermediaire}
\mu \left( A_{\alpha}(S) \right) \leqslant 9 \eta
+ \mu \left( \left\{ \exists s \in g_S(\omega) \cap \bigcup_{k = [\eta S]}^{S}
\left( L_k \cup R_k \right) :  \ell(s) \leqslant S^{\alpha} \right\} \cap \left\{
\forall k \geqslant [\eta S],  Y_k \geqslant [\delta \sqrt{S}] \right\}
\right).
\end{equation}
Let us fix a collection $(y_k)_{[\eta S] \leqslant k \leqslant S}$ such that
$y_k \geqslant [\delta \sqrt{S}]$ for every $k$. Theorem \ref{descmu} gives
\begin{align}
\label{conditionel}
\mu & \left( \exists s \in g_S(\omega) \cap \bigcup_{k = [\eta S]}^{S}
\left( L_k \cup R_k \right) \, : \, \ell(s) \leqslant S^{\alpha} \middle|
Y_k = y_k, \, [\eta S] \leqslant k \leqslant S \right) \notag \\
& \leqslant
2 \sum_{k = [\eta S]}^{S}
\widehat{\rho}^{(y_k)} \left( \vphantom{\sqrt{S}} \exists s \in g_{S-k} (\omega) :
  \, \ell(s) \leqslant S^{\alpha} \right)
 = 2 \sum_{k=0}^{S - [\eta S]}
\widehat{\rho}^{(y_{S-k})} \left(\vphantom{\sqrt{S}} \exists s \in g_{k} (\omega) :
  \, \ell(s) \leqslant S^{\alpha} \right).
\end{align}
If $0 \leqslant k \leqslant S - [\eta S]$, one has:
\begin{equation*}
\begin{split}
\widehat{\rho}^{(y_{S-k})} \left(\vphantom{\sqrt{S}} \exists s \in g_{k} (\omega) :
  \, \ell(s) \leqslant S^{\alpha} \right) &  \leqslant
\widehat{\rho}^{(y_{S-k})} \left(\vphantom{\sqrt{S}}
\inf_{s \in \omega} \ell(s)  \leqslant S^{\alpha} \right)  \\
& = \frac{1}{w_{y_{S-k}}} \sum_{\substack{\omega \in \mathbb{T}^{(y_{S-k})}\\ \inf_{s \in \omega}
  \ell(s)  \leqslant S^{\alpha}}} 12^{- |\omega|} \\
& = \frac{1}{w_{y_{S-k}}} \sum_{\substack{\omega \in
  \mathbf{T}^{(y_{S-k})}\\ 0 < \inf_{s \in \omega}
  \ell(s)  \leqslant S^{\alpha} }} 12^{- |\omega|}\\
& = \frac{1}{w_{y_{S-k}}} \rho^{(y_{S-k})} \left(\vphantom{\sqrt{S}}
 0 < \inf_{s \in \omega} \ell(s) \leqslant S^{\alpha} \right) .\\
\end{split}
\end{equation*}
But
\[
\rho^{(y_{S-k})} \left(\vphantom{\sqrt{S}}
\inf_{s \in \omega} \ell(s) > 0 \right) = w_{y_{S-k}} 
\]
and
\begin{equation*}
\rho^{(y_{S-k})} \left(\vphantom{\sqrt{S}}
\inf_{s \in \omega} \ell(s) > S^{\alpha} \right) 
= \rho^{(y_{S-k} - [S^{\alpha}])} \left(\vphantom{\sqrt{S}} \inf_{s
    \in \omega} \ell(s) > 0 \right)
= w_{y_{S-k} - [S^{\alpha}]}.
\end{equation*}
We thus have
\[
\widehat{\rho}^{(y_{S-k})} \left(\vphantom{\sqrt{S}} \exists s \in g_{k} (\omega) :
  \, \ell(s) \leqslant S^{\alpha} \right)  \leqslant
\frac{1}{w_{y_{S-k}}} \left( w_{y_{S-k}} - w_{y_{S-k} -  [S^{\alpha}]}
\right) = 1 - \frac{w_{y_{S-k} - [S^{\alpha}]}}{w_{y_{S-k}}}.
\]
Using our assumption $y_{S-k} \geqslant [ \delta \sqrt{S}]$, a Taylor expansion gives
\[
1 - \frac{w_{y_{S-k} - [S^{\alpha}]}}{w_{y_{S-k}}} =
4 S^{\alpha} y_{S-k}^{-3} + o\left( S^{\alpha} y_{S-k}^{-3} \right)
 \leqslant \frac{4}{\delta^3} S^{\alpha - \frac{3}{2}} +  o\left( S^{\alpha - \frac{3}{2}}  \right),\]
and the right-hand side of \eqref{conditionel} is smaller than
$\frac{8}{\delta^3} S^{\alpha -1/2} + o \left(S^{\alpha - 1/2} \right)$.
From \eqref{Aalphaintermediaire} we now get
\[
\mu \left( A_{\alpha}(S) \right) \leqslant 9 \eta +
\frac{8}{\delta^3} S^{\alpha - \frac{1}{2}} +  o\left( S^{\alpha - \frac{1}{2}} \right).
\]
Hence $\mu(A_{\alpha}(S)) < 10 \eta$ as soon as $S$ is large enough. This completes the proof.
\end{proof}

\begin{proposition}
\label{nseuil}
For every sufficiently large integer $S$, there exists an integer $N(S)$ such that, for every $n>N(S)$, one has
\[
\mu_n(\Omega_S) \leqslant \varepsilon.
\]
\end{proposition}
\begin{proof}
In this proof $\alpha \in \left]\frac{1}{3}, \frac{1}{2} \right[$ is fixed. For $S > 0$, Lemma \ref{size} gives $K_{\varepsilon}(S) > 0$ and $N_1(S) > 0$ such that if $n > N_1(S)$ then $\mu \left( \omega : \, |B_{\overline{\mathbb{T}},S} (\omega)| > K_{\varepsilon}(S)  \right) < \varepsilon $. Let us also recall that the number of well-labelled trees with height $S$ and size smaller than $K_{\varepsilon}(S)$ is denoted by $M_{\varepsilon}(S)$. Lemma \ref{estimlabels} shows that, for $S$ large enough, there exists $N_2(S)$ such that $\mu_n(A_{\alpha}(S)) < \varepsilon$ for every $n > N_2(S)$. Therefore, for $S$ large enough and for $n > N_1(S) \vee N_2(S)$ one has
\begin{align}
\mu_n (\Omega_S)  & = \sum_{\substack{\omega^{\star} \notin A_{\alpha}(S)\\
|\omega^{\star}| \leqslant K_{\varepsilon}(S), \, h(\omega^{\star}) = S}} \mu_n \left( \{\omega 
  : \, B_{\overline{\mathbb{T}},S} (\omega) = \omega^{\star} \} \cap
\Omega_S \right) \notag\\
& \qquad + \mu_n (A_{\alpha}(S)) +  \mu_n
\left( \omega: \, |B_{\overline{\mathbb{T}},S} (\omega)| > K_{\varepsilon}(S)
\right) \notag \\
& \leqslant 2 \varepsilon + \sum_{\substack{\omega^{\star} \notin A_{\alpha}(S)\\
|\omega^{\star}| \leqslant K_{\varepsilon}(S), \, h(\omega^{\star}) = S}} \mu_n \left( \{\omega 
  : \, B_{\overline{\mathbb{T}},S} (\omega) = \omega^{\star} \} \cap
\Omega_S  \right).
\label{munosr}
\end{align}

Fix a tree $\omega^{\star} \notin A_{\alpha}(S)$ with height $S$ and size smaller than $K_{\varepsilon}(S)$. We assume that $S$ is large enough so that $S^{\alpha} > R + 1$. We denote by $k$ the number of vertices of $\omega^{\star}$ at generation $S$ and by $l_1, \ldots , l_k$ the labels of these vertices. By considering the subtrees of $\omega$ originating from vertices at generation $S$, one obtains:
\begin{equation}
\label{omegastar}
\mu_n \left( \{\omega : \, B_{\overline{\mathbb{T}},S}(\omega) = \omega^{\star} \} \cap
\Omega_S \right) \leqslant
\frac{1}{D_n} \, \sum_{n_1 + \cdots + n_k = n - |\omega^{\star}|} \, \sum_{i = 1}^k
D_{n_i}^{(l_i)}(R) \prod_{j \neq i} D_{n_j}^{(l_j)}
\end{equation}
where $D_n^{(l)}(R)$ is the number of trees in $\mathbb{T}_n^{(l)}$
with at least one vertex with label less than or equal to $R+1$ (compare \eqref{omegastar} with formula \eqref{boulen}).
Since $\omega \notin A_{\alpha}(S)$, we have
$l_i > S ^{\alpha} > R+1$ and thus $D_{n_i}^{(l_i)}(R) = D_{n_i}^{(l_i)} - D_{n_i}^{(l_i - R - 1)}$
for $i = 1, \ldots , k$. The bound \eqref{omegastar} then gives
\begin{align*}
\mu_n  & \left( \{\omega : \, B_{\overline{\mathbb{T}},S}(\omega) = \omega^{\star} \} \cap
\Omega_S \right) \\
& \leqslant
\frac{1}{D_n} \, \sum_{n_1 + \cdots + n_k = n - |\omega^{\star}|} \, \sum_{i = 1}^k
(D_{n_i}^{(l_i)} - D_{n_i}^{(l_i - R - 1)}) \prod_{j \neq i} D_{n_j}^{(l_j)}\\
& =
\frac{k}{D_n} \, \sum_{n_1 + \cdots + n_k = n - |\omega^{\star}|} \prod_{j=1}^k
D_{n_j}^{(l_j)}
- \frac{1}{D_n} \,  \sum_{i = 1}^k \, \sum_{n_1 + \cdots + n_k = n -
  |\omega^{\star}|} \, D_{n_i}^{(l_i - R - 1)} \prod_{j \neq i} D_{n_j}^{(l_j)}.
\end{align*}

By Theorem \ref{cvarbres}, $\mu_n \left(\omega : B_{\overline{\mathbb{T}},S} (\omega) =
\omega^{\star} \right) \rightarrow \mu \left(\omega : B_{\overline{\mathbb{T}},S} (\omega) =
\omega^{\star} \right)$  as $n \rightarrow \infty$. Using this convergence and identities \eqref{boulen} and \eqref{boule}, we get the existence of an integer $N(\omega^{\star},S)$ such that for $n > N(\omega^{\star},S)$ one has
\begin{equation*}
\frac{1}{D_n} \,  \sum_{n_1 + \cdots + n_k = n - |\omega^{\star}|} \prod_{j=1}^k
D_{n_j}^{(l_j)}
\leqslant
12^{- |\omega^{\star}|} \, \sum_{t = 1}^k d_{l_t} \prod_{s \neq t} w_{l_s} \, +
\, \frac{\varepsilon}{K_{\varepsilon}(S)M_{\varepsilon}(S)}
\end{equation*}
and for $i = 1 ,\ldots, k$
\begin{equation*}
\begin{split}
 \frac{1}{D_n}\, & \sum_{n_1 + \cdots + n_k = n - |\omega^{\star}|} \,
D_{n_i}^{(l_i - R -1)} \prod_{j \neq i} D_{n_j}^{(l_j)} \\
& \geqslant
12^{- |\omega^{\star}|} \Big( d_{l_i - R -1} \prod_{j \neq i} w_{l_j} + \sum_{t \neq
  i} d_{l_t} w_{l_i - R -1} \prod_{j \neq t, i} w_{l_j} \Big) -
\frac{\varepsilon}{K_{\varepsilon}(S)M_{\varepsilon}(S)}.
\end{split}
\end{equation*}
We now have for every $n > N(\omega^{\star},S)$:
\begin{align}
\mu_n & \left( \{\omega : \, B_{\overline{\mathbb{T}},S}(\omega) = \omega^{\star} \} \cap
\Omega_S \right) \notag \\
 & \leqslant \frac{2\varepsilon}{M_{\varepsilon}(S)}
+ k 12^{- |\omega^{\star}|} \, \sum_{t = 1}^k d_{l_t} \prod_{s \neq t} w_{l_s}
- 12^{- |\omega^{\star}|} \sum_{i=1}^k \Big( d_{l_i - R -1} \prod_{j \neq i} w_{l_j} + \sum_{t \neq
  i} d_{l_t} w_{l_i - R -1} \prod_{j \neq t, i} w_{l_j} \Big) \notag \\
 & = \frac{2\varepsilon}{M_{\varepsilon}(S)}
+ 12^{- |\omega^{\star}|} \, \sum_{t = 1}^k (d_{l_t} - d_{l_t -R-1}) \prod_{s \neq
  t} w_{l_s}
+ 12^{- |\omega^{\star}|} \, \sum_{t = 1}^k d_{l_t}\left( \sum_{i \neq
    t} (w_{l_i} - w_{l_i - R-1}) \prod_{s \neq t,i} w_{l_s} \right).\label{muncapomega}
\end{align}
Define
\begin{align*}
d(\omega^{\star}) & = \max_{i= 1 \ldots k} \left( 1 - \frac{d_{l_i - R-1}}{d_{l_i}}
\right),\\
w(\omega^{\star}) & =  \max_{i= 1 \ldots k} \left( 1 - \frac{w_{l_i - R-1}}{w_{l_i}}
\right).
\end{align*}
From the bound \eqref{muncapomega}, we get
\begin{align}
\mu_n & \left( \{\omega : \, B_{\overline{\mathbb{T}},S} (\omega) = \omega^{\star} \} \cap
\Omega_S \right) \notag \\
 & \leqslant \frac{2\varepsilon}{M_{\varepsilon}(S)}
+ d(\omega^{\star}) \, 12^{- |\omega^{\star}|} \, \sum_{t = 1}^k d_{l_t} \prod_{s \neq
  t} w_{l_s}
+ k w(\omega^{\star}) \, 12^{- |\omega^{\star}|} \, \sum_{t = 1}^k d_{l_t} \prod_{s \neq
  t} w_{l_s} \notag \\
& = \frac{2\varepsilon}{M_{\varepsilon}(S)}
+ \left( d(\omega^{\star}) + k . w(\omega^{\star}) \right) \mu \left( \omega :
B_{\overline{\mathbb{T}},S} (\omega) = \omega^{\star} \right)
\label{omegastarfinal}
\end{align}
where we used \eqref{boule} in the last equality.

\bigskip

Let us now define $N^{\star}(S) = \max_{|\omega^{\star}| \leqslant K_{\varepsilon}(S)} N(\omega^{\star},S) \vee
N_1(S) \vee N(S)$. For $S$ large enough and for $n > N^{\star}(S)$ we obtain using (\ref{munosr}) and (\ref{omegastarfinal}):
\begin{equation*}
\mu_n(\Omega_S) \leqslant 4 \varepsilon
+ \sum_{\substack{\omega^{\star} \notin A_{\alpha}(S)\\ |\omega^{\star}| \leqslant
  K_{\varepsilon}(S), \, h(\omega^{\star}) = S}}
\left( d(\omega^{\star}) + |g_S(\omega^{\star})| . w(\omega^{\star}) \right) \mu \left(\omega :
B_{\overline{\mathbb{T}},S} (\omega) = \omega^{\star} \right).
\end{equation*}

A Taylor expansion gives $w(\omega^{\star}) \leqslant 4 (5R  + 2) S^{-3 \alpha} + o\left(S^{-3 \alpha} \right)$ where the remainder is uniform over $\omega^{\star} \notin A_{\alpha}(S)$. In addition, $\sup_{\omega^{\star} \notin A_{\alpha}(S)} d(\omega^{\star}) \rightarrow 0$ as $S \rightarrow \infty$. This allows us to find $S^{\star}$ such that for $S > S^{\star}$ and $n> N^{\star}(S)$:
\begin{align}
\mu_n(\Omega_S) & \leqslant 4 \varepsilon +
\sum_{\substack{\omega^{\star} \notin A_{\alpha}(S)\\  |\omega^{\star}| \leqslant K_{\varepsilon}(S), \, h(\omega^{\star}) = S}}
\left( \varepsilon + |g_S(\omega^{\star})| . 4(5R+2) S^{-3 \alpha} \right) \mu \left(\omega :
B_{\overline{\mathbb{T}},S} (\omega) = \omega^{\star} \right) \notag \\
& \leqslant 5 \varepsilon +
4(5R+2) S^{-3 \alpha} \sum_{\substack{\omega^{\star} \notin A_{\alpha}(S)\\ |\omega^{\star}| \leqslant
  K_{\varepsilon}(S), \, h(\omega^{\star})=S}} |g_S(\omega^{\star})|  \, \mu \left(\omega :
B_{\overline{\mathbb{T}},S} (\omega) = \omega^{\star} \right) \notag \\
& \leqslant 5 \varepsilon +
4(5R+2) S^{-3 \alpha} \mathbb{E}_{\mu}
\left[|g_S(\omega)| \right]. \label{munomegas}
\end{align}

The description of $\mu$ given in Theorem \ref{descmu} allows us to estimate $\mathbb{E}_{\mu} \left[|g_S(\omega)| \right]$. Indeed we have for every integer $H > 0$ and $k \geqslant 1$
\[
\mathbb{E}_{\widehat{\rho}^{(k)}} \left[ |g_H(\omega)| \right] \leqslant
\frac{1}{w_k} \mathbb{E}_{\rho^{(k)}} \left[ |g_H(\omega)| \right] =
\frac{2}{w_k} \mathbb{E}_{GW(1/2)} \left[ |g_H(\omega)| \right] =
\frac{2}{w_k} \leqslant 2.
\]
It follows that
\[
\mathbb{E}_{\mu} \left[ |g_S(\omega)| \right] \leqslant 4S +1.
\]
Recalling that $\alpha > \frac{1}{3}$, we get that for every $S$ large enough
and for $n > N^{\star}(S)$,
\[
\mu_n(\omega_S) \leqslant 6 \varepsilon.
\]
This completes the proof.
\end{proof}

\subsection{Proof of the main result}

In this section we fix $Q^{\star} \in \overline{\mathbf{Q}}$ and $R >0$. As in the previous section, we write $\Omega_S = \Omega_S(R)$ to simplify notation. From the remarks following Theorem \ref{equality}, the proof reduces to verifying the convergence
\begin{equation}
\label{conv}
\mu_n \Big( \omega :
B_{\overline{\mathbf{Q}},R} (\Phi(\omega)) = B_{\overline{\mathbf{Q}},R} (Q^{\star}) \Big)
\underset{n \to \infty}{\longrightarrow}
\mu \Big( \omega :
B_{\overline{\mathbf{Q}},R} (\Phi(\omega)) = B_{\overline{\mathbf{Q}},R} (Q^{\star}) \Big).
\end{equation}

First of all, we need to reformulate the problem in terms of trees. Since $Q^{\star} \in \overline{\mathbf{Q}}$, we know that there
exists a finite quadrangulation $Q_0 \in \mathbf{Q}$ such that  $d_{\mathbf{Q}} (Q_0,Q^{\star}) <
\frac{1}{R+1}$ and therefore $B_{\overline{\mathbf{Q}},R} (Q_0) = B_{\overline{\mathbf{Q}},R} (Q^{\star})$.
Then there exists $\omega_0 \in \mathbb{T}$ such that $\Phi(\omega_0) = Q_0$. The convergence (\ref{conv}) can now be restated as
\begin{equation}
\label{convbis}
\mu_n \Big( \omega :
B_{\overline{\mathbf{Q}},R} (\Phi(\omega)) = B_{\overline{\mathbf{Q}},R} (\Phi(\omega_0)) \Big)
\underset{n \to \infty}{\longrightarrow} \mu \Big( \omega :
B_{\overline{\mathbf{Q}},R} (\Phi(\omega)) = B_{\overline{\mathbf{Q}},R} (\Phi(\omega_0)) \Big).
\end{equation}
We fix $\varepsilon > 0$ in the remaining part of this proof.

\bigskip

We need to characterize the trees $\omega$ for which $\Phi(\omega)$ has the same ball of radius $R$ as $\Phi(\omega_0)$. As we have already mentioned at the end of Section \ref{sec:quadrangulations}, the main difficulty comes from the fact that two trees that are very similar in $\overline{\mathbb{T}}$ can give very different quadrangulations if they have vertices with small labels in high generations. We can remedy this problem thanks to Proposition \ref{nseuil}.

Note that $\omega_0$ is a finite tree. Let $S_0$ denote the height of $\omega_0$. According to Lemma \ref{seuil} and Proposition \ref{nseuil} we can choose $S_1 > S_0$ such that if $S \geqslant S_1$ and $n \geqslant N(S)$ then $\mu (\Omega_S) < \varepsilon$ and $\mu_n (\Omega_S) < \varepsilon$.

Let $S > S_1$ and let $(\omega_i)_{i \in I}$ be the collection of trees given by Corollary \ref{pbcle}, such that, for every $\omega \in \mathscr{S} \cap \Omega_S^c$, the equality $B_{\overline{\mathbf{Q}},R} (\Phi(\omega)) =
B_{\overline{\mathbf{Q}},R} (\Phi(\omega_0))$ holds if and only if there exists $i \in I$ such that
$B_{\overline{\mathbb{T}},S} (\omega) = B_{\overline{\mathbb{T}},S} (\omega_i)$. If $A
\vartriangle B$ denotes the symmetric difference between two sets $A$ and $B$, we have
\begin{equation*}
\begin{split}
\mu & \left( \left\{ \omega \in \overline{\mathbb{T}} : \,
B_{\overline{\mathbf{Q}},R} (\Phi(\omega)) = B_{\overline{\mathbf{Q}},R} (\Phi(\omega_0)) \right\}
\vartriangle
\bigcup_{i\in I} \left\{ \omega \in \overline{\mathbb{T}} : \, 
B_{\overline{\mathbb{T}},S} (\omega) = B_{\overline{\mathbb{T}},S} (\omega_i) \right\}
\right)\\
& \leqslant \mu \left( \Omega_S \right) < \varepsilon.
\end{split}
\end{equation*}
We deduce from this last bound that
\begin{equation*}
\begin{split}
\Bigg| & \mu_n \Big( \omega : \,
B_{\overline{\mathbb{Q}},R} (\Phi(\omega)) = B_{\overline{\mathbf{Q}},R} (\Phi(\omega_0)) \Big)
- \mu \Big( \omega : \,
B_{\overline{\mathbf{Q}},R} (\Phi(\omega)) = B_{\overline{\mathbf{Q}},R} (\Phi(\omega_0)) 
\Big) \Bigg|\\
& \leqslant  \Bigg| \mu \Big(
\bigcup_{i\in I} \left\{ \omega : \,
B_{\overline{\mathbb{T}},S} (\omega) = B_{\overline{\mathbb{T}},S} (\omega_i) \right\} \Big)
- \mu_n \Big(\omega : \,
B_{\overline{\mathbf{Q}},R} (\Phi(\omega)) = B_{\overline{\mathbf{Q}},R}
(\Phi(\omega_0)) \Big) \Bigg|\\
 &  \qquad + \Bigg| \mu \Big(
\bigcup_{i\in I} \left\{ \omega : \,
B_{\overline{\mathbb{T}},S} (\omega) = B_{\overline{\mathbb{T}},S} (\omega_i) \right\} \Big)
- \mu \Big( \omega : \,
B_{\overline{\mathbf{Q}},R} ( \Phi(\omega)) = B_{\overline{\mathbf{Q}},R}
(\Phi(\omega_0)) \Big) \Bigg|\\
 & \leqslant \Bigg| \mu \Big(
\bigcup_{i\in I} \left\{ \omega : \,
B_{\overline{\mathbb{T}},S} (\omega) = B_{\overline{\mathbb{T}},S} (\omega_i) \right\} \Big)
- \mu_n \Big( \omega : \,
B_{\overline{\mathbf{Q}},R} (\Phi(\omega)) = B_{\overline{\mathbf{Q}},R} (\Phi(\omega_0)) 
\Big) \Bigg| + \varepsilon .
\end{split}
\end{equation*}
The set $\bigcup_{i\in I} \left\{ \omega \in \overline{\mathbb{T}} : \, 
B_{\overline{\mathbb{T}},S} (\omega) = B_{\overline{\mathbb{T}},S} (\omega_i) \right\}$
is both open and closed in $\overline{\mathbb{T}}$, and thus
\begin{equation*}
\mu_n \left( \bigcup_{i\in I} \left\{
\omega : \, B_{\overline{\mathbb{T}},S} (\omega) 
= B_{\overline{\mathbb{T}},S}(\omega_i) \right\} \right)
\underset{n \to \infty}{\longrightarrow}
\mu \left( \bigcup_{i\in I} \left\{
\omega : \, B_{\overline{\mathbb{T}},S} (\omega) 
= B_{\overline{\mathbb{T}},S} (\omega_i) \right\} \right).
\end{equation*}
Therefore there exists $N'(S) > 0$ such that for $n > N'(S)$:
\begin{equation*}
\begin{split}
\Bigg| & \mu_n \Big( \omega : \,
B_{\overline{\mathbf{Q}},R} (\Phi(\omega)) = B_{\overline{\mathbf{Q}},R} (\Phi(\omega_0)) \Big)
- \mu \Big( \omega : \,
B_{\overline{\mathbf{Q}},R} (\Phi(\omega)) = B_{\overline{\mathbf{Q}},R} (\Phi(\omega_0)) \Big) \Bigg|\\
 & \leqslant \Bigg| \mu_n \Big( \bigcup_{i\in I} \left\{ \omega : \,
B_{\overline{\mathbb{T}},S} (\omega) = B_{\overline{\mathbb{T}},S} (\omega_i) \right\} \Big)
- \mu_n \Big( \omega : \,
B_{\overline{\mathbf{Q}},R} (\Phi(\omega)) = B_{\overline{\mathbf{Q}},R} (\Phi(\omega_0)) \Big) \Bigg|
+ 2 \varepsilon \\
 & = \Bigg| \mu_n \Big( \bigcup_{i\in I} \left\{ \omega : \,
B_{\overline{\mathbb{T}},S} (\omega) = B_{\overline{\mathbb{T}},S} (\omega_i) \right\} \cap
\Omega_S \Big)
- \mu_n \Big( \left\{ \omega : \,
B_{\overline{\mathbf{Q}},R} (\Phi(\omega)) = B_{\overline{\mathbf{Q}},R} (\Phi(\omega_0)) 
\right\} \cap \Omega_S \Big) \Bigg| \\
& \qquad + 2 \varepsilon
\end{split}
\end{equation*}
by the choice of the collection $(\omega_i)_{i \in I}$.

\bigskip

We also know that $\mu_n(\Omega_S) < \varepsilon$ for $n > N(S)$, and it follows that, for $n > N(S) \vee N'(S)$,
\begin{equation*}
\Bigg| \mu_n
\Big( \omega : \,
B_{\overline{\mathbf{Q}},R} (\Phi(\omega)) = B_{\overline{\mathbf{Q}},R} (\Phi(\omega_0)) \Big)
- \mu
\Big( \omega : \,
B_{\overline{\mathbf{Q}},R} (\Phi(\omega)) = B_{\overline{\mathbf{Q}},R} (\Phi(\omega_0)) \Big) \Bigg|
\leqslant 3 \varepsilon.
\end{equation*}
This completes the proof of Theorem \ref{equality}. \qed

\bigskip

\noindent {\bf Acknowledgements.} The author would like to thank Jean-Fran\c{c}ois Le Gall for many helpful discussions about this work.

\addcontentsline{toc}{section}{References}
\bibliographystyle{abbrv}
\def\cprime{$'$}

\end{document}